\documentclass[a4paper,10pt]{amsart}
\usepackage[matrix,arrow,tips,curve]{xy}
\usepackage[english]{babel}
\usepackage[leqno]{amsmath}
\usepackage{amssymb,amsthm}
\usepackage{amscd}
\usepackage{enumerate}

\input{xy}
\xyoption{all}
\newtheorem{thm}{Theorem}[section]
\newtheorem{lemma}[thm]{Lemma}

\newtheorem{prop}[thm]{Proposition}

\newtheorem{cor}[thm]{Corollary}

\newtheorem{fact}[thm]{Fact}

\theoremstyle{remark}
\newtheorem{remark}[thm]{Remark}
\theoremstyle{definition}
\newtheorem{defi}[thm]{Definition}
\newtheorem{nota}[thm]{}

\numberwithin{equation}{section}

\newtheorem{example}[thm]{Example}

\newcommand{\la}{\longrightarrow}

\newcommand{\ov}{\overline}

\newcommand{\pr}[1]{\mathbb{P}^{#1}}

\newcommand{\im}{\operatorname{Im}}
\newcommand{\codim}{\operatorname{codim}}

\newcommand{\Pic}{\operatorname{Pic}}

\newcommand{\mdeg}{\operatorname{{\underline{de}g}}}

\def\O{\mathcal O}

\def\X{\mathcal X}

\newcommand{\Z}{\mathbb{Z}}

\newcommand{\ph}{\varphi}

\def\md{\underline{d}}
\def\mci{\underline{c}_i}

\def\PX{P^d_X}
\def\PXb{\overline{P^d_X}}

\def\pf{P^d_f}

\def\pfb{\overline{P^d_f}}

\newcommand{\sing}{X_{\text{sing}}}
\newcommand{\sep}{X_{\text{sep}}}

\def\nf{N^d_f}
\def\LX{\Lambda _X}

\def\DX{\Delta _X}
\def\Ddg{D_g^d}

\def\Mgbst{\overline{\mathcal{M}}_g}

\newcommand{\XS}{X_S^{\nu}}

\newcommand{\eX}{X_{\text{exc}}}

\newcommand{\hXS}{{\widehat{X}}_S}
\newcommand{\nXS}{{X}^{\nu}_S}

\newcommand{\mdn}{{\md}^{\nu}}
\newcommand{\hL}{{\hat{L}}}

\newcommand{\Mgb}{\ov{M_g}}
\newcommand{\Pdgb}{\ov{P_{d,g}}}

\begin{document}

\title{Compactified Jacobians of N\'eron type}
\author{Lucia Caporaso}
\address{Dipartimento di Matematica,
Universit\`a Roma Tre,
Largo S. Leonardo Murialdo 1,
00146 Roma (Italy)}
\email{caporaso@mat.uniroma3.it}
\maketitle

 \begin{abstract}
We characterize   stable curves  $X$   whose compactified degree-$d$  Jacobian
is of N\'eron type.
This means the following: for any one-parameter regular smoothing of $X$, 
the special fiber of the N\'eron model of the   Jacobian of the generic fiber
is isomorphic to a dense open subset of the degree-$d$ compactified Jacobian.
It is well known that compactified Jacobians of N\'eron type have the best modular properties, and
that they are endowed with a mapping property useful for applications.
\end{abstract}

 \section{Introduction and preliminaries}
 
 Let $X$ be a stable curve and $f:\X\to B$ a one-parameter smoothing of $X$ with
 $\X$ a nonsingular surface;  $X$ is thus identified with the special fiber of $f$ and all other fibers are smooth curves.
 Let $N^d_f\to B$ be the N\'eron model of the degree-$d$ Jacobian of the generic fiber of $f$;
its existence was proved by A. N\'eron in \cite{neron}, and its connection  with the Picard functor
was  established by M. Raynaud in \cite{raynaud}. 
 So, $N^d_f\to B$  is a smooth and separated morphism, 
 whose generic fiber is the degree-$d$ Jacobian of the generic fiber of $f$; the special fiber, denoted $N^d_X$, is isomorphic to a disjoint union of copies of the generalized Jacobian of $X$. 
$N^d_f\to B$   can be interpreted as the maximal separated quotient of the degree-$d$ Picard scheme $\Pic^d_f\to B$.
In particular, if $\Pic^d_f\to B$ is separated, 
 which happens if and only if $X$ is irreducible,  then $N^d_f=\Pic^d_f$
 (we refer to \cite{raynaud}, \cite{BLR} or \cite{artin} for details).
 
 The N\'eron model 
 has a universal property,
  the N\'eron Mapping Property 
 (\cite[ def. 1]{BLR}),
which determines it uniquely. 
Moreover, as $d$ varies in $\Z$, the special fibers, $N^d_X$, of $N^d_f\to B$ are all isomorphic.
 
By contrast, the compactified degree-$d$ Jacobian 
 of a reducible curve $X$, denoted  $\PXb$,
 has a structure which varies with $d$. For example,     the number of irreducible components,
 and the modular properties, 
 depend on $d$; see Section~\ref{sec2} for details and references.

$\PXb$  will be called  of {\it N\'eron type}
 if its smooth locus is isomorphic to $N^d_X$.
 Compactified Jacobians of N\'eron type have the best modular properties.
Moreover they inherit a mapping property from the universal property of the N\'eron model  which provides a very
useful tool; see for example  \cite{CE} for applications to Abel maps.  

The purpose of this paper is to classify, for every $d$, those  stable curves 
$X\in \Mgb$ such that   $\PXb$ is on N\'eron type.
The question is interesting if $g\geq 2$, for otherwise $\PXb$ is always irreducible, and hence of N\'eron type.

Before stating our main result, we need a few words about compactified Jacobians.
$\PXb$  parametrizes certain line bundles on quasistable curves having $X$ as stabilization.
These are the so-called ``balanced"  line bundles;
among balanced line bundles there are some distinguished ones, called ``strictly balanced", which have better modular properties.
In fact, to every balanced line bundle there corresponds a unique point in 
$\PXb$, but different balanced line bundles may determine the same point.
On the other hand every point of $\PXb$ corresponds to a unique
  class of strictly balanced line bundles.

  The curve $X$ is called {\it $d$-general} if every balanced line bundle of degree $d$ is strictly balanced. 
 This is equivalent to the fact that $\PXb$ is a geometric GIT-quotient.
 The property of being $d$-general depends only on the weighted dual graph of $X$,
 and the locus of $d$-general curves in $\Mgb$ has been precisely described by M. Melo in \cite{melo}.

Now, the degree-$d$ compactified Jacobian of a $d$-general curve is of N\'eron type, by \cite[Thm. 6.1]{cner}. 
 But, as we will prove, the converse does not hold.
 
More precisely,     a stable curve $X$ is called {\it weakly $d$-general} if a    curve obtained by smoothing every separating node of $X$, and maintaining all the non separating nodes, is $d$-general; see Definition~\ref{wdef}.  

Our main result, Theorem~\ref{main},  states that $\PXb$ is of N\'eron type if and only if $X$ is weakly $d$-general. The proof combines standard techniques with combinatorial methods originating from  recent  work with F. Viviani.

The locus of weakly $d$-general curves in $\Mgb$ is precisely described in section~\ref{locus};
its complement turns out to have codimension at least $2$.

I am grateful to  M. Melo and F. Viviani for their precious comments.
\begin{nota}{\it Notations and conventions.}
\label{not}
\begin{enumerate}
\item
We work over an algebraically closed field $k$.
The word ``curve" means projective scheme of pure dimension one.
The genus of a curve will be the arithmetic genus, unless otherwise specified.

\item
\label{notZ}
By  $X$ we will always denote  a nodal  curve of genus $g$. 
For any subcurve $Z\subset X$   we denote by $g_Z$   its arithmetic genus, by
$Z^c:=\overline{X\smallsetminus Z}$ and by
$\delta _Z:=\#Z\cap Z^c$. 
We set
$ 
w_Z:=\deg_Z\omega_X=2g_Z-2+\delta_Z.
$ 
\item
A node $n$ of a connected curve $X$ is called {\it separating} if $X\smallsetminus \{n\}$
is not connected. The set of all separating nodes of $X$ is denoted by $\sep$ and the set of all nodes of $X$ by $\sing$.

\item
A nodal curve $X$  of genus $g\geq 2$ is called {\it stable} if it is connected and if
every   component $E\subset X$ such that $E\cong  \pr{1}$ satisfies $\delta_E\geq 3$.
$X$ is called {\it quasistable} if it is connected,
 if
every  $E\subset X$ such that $E\cong \pr{1}$ satisfies $\delta_E\geq 2$,
and if two   exceptional components never intersect, where an exceptional component is defined as an $E\cong \pr{1}$ such that $\delta_E=2$.
We denote by $\eX$ the union of   the exceptional components of $X$.

\item
\label{notS}
Let $S\subset \sing$
we denote by $\nu_S:\nXS\to X$ the normalization of $X$ at $S$, and by $\hXS$ the quasistable curve obtained by ``blowing-up" all the nodes in $S$, so that there is a natural surjective map
$$ 
\hXS=\bigcup_{i=1}^{\#S}E_i\cup\nXS\la X
$$ 
restricting to $\nu_S$ on $\XS$ and
contracting all the exceptional components $E_i$ of $\hXS$.
 $\hXS$ is also called a {\it quasistable curve of} $X$.
\item
Let $C_1,\ldots, C_{\gamma}$ be the irreducible components of $X$. Every line bundle $L\in \Pic X$ has a multidegree $\mdeg L=(deg_{C_1}L,\ldots,deg_{C_{\gamma}})\in \Z^{\gamma}$.
Let $\md=(d_1,\ldots, d_{\gamma})\in \Z^{\gamma}$, we set
$|\md|=\sum_1^{\gamma}d_i$; for any
subcurve  $Z\subset X$    we abuse notation slightly and denote
$$
\md_Z:=\sum_{C_i\subset Z}d_i.
$$
\end{enumerate}
\end{nota}

\begin{nota}{\it Compactified Jacobians of N\'eron type.}
\label{neron}
Let  $X$ be any nodal connected curve and $f:\X \to B$ a one-parameter {\it regular smoothing} for $X$,
i.e. $B$ is a smooth connected one-dimensional scheme 
 with a marked point $b_0\in B$, $\X$ is a regular surface, and $f$ is a projective morphism  whose fiber over $b_0$ is $X$ and whose remaining fibers are smooth curves.
We set $U:=B\smallsetminus\{b_0\}$ and let $f_U:\X_U\to U$ be the family of smooth curves obtained by restricting $f$ to $U$.
Consider the relative degree $d$ Picard scheme
over $U$, denoted
$ 
\Pic^d_{f_U}\to U.
$ 
Its N\'eron model over $B$
will be denoted 
 by
\begin{equation}
\label{neron}
\nf:=N(\Pic^d_{f_U})\la B.
\end{equation}
The fiber of $\nf$ over $b_0$ will be denoted $N_X^d$;
$N_X^d$ is isomorphic to a finite number of copies of the generalized Jacobian of $X$. The number of copies is independent of $d$; to compute it we
 introduce the so-called ``degree class group".
 
  Let $\gamma$ be the number of irreducible components of $X$.
For every component $C_i$ of $X$ set 
$
k_{i,j}:= \#(C_i\cap C_j)
$  if $j \neq i$,
and 
$
k_{i,i}=-\#(C_i\cap \overline{C\setminus C_i})
$
so that the  matrix
$
(k_{i,j})
$
is    symmetric matrix.
Notice that for every regular smoothing $f:\X\to B$ of $X$ as above, we have
$
\deg _{C_j}\O_\X(C_i)=k_{i,j}.
$
Hence this matrix is also related to $f$, although it does not depend on the choice of $f$
(as long as $\X$ is regular).

We have $\sum _{j=1}^{\gamma} k_{i,j}=0$ for every $i$.
Now, for every $i=1,\ldots ,\gamma$ set 
$
\mci :=(k_{1,i},\ldots ,k_{\gamma ,i})\in \Z ^{\gamma}
$\  
and
$ 
{\bf Z}:=\{ \md \in \Z ^{\gamma}: \  |\md | = 0\}
$ 
so that $\mci\in {\bf Z}$. We can now define the  sublattice  
$ 
\LX :=< \underline{c}_1, \ldots ,\underline{c}_{\gamma} >\subset {\bf Z}.
$

The {\it degree class group} of $X$ is the  group
$\DX := {\bf Z}/\LX$. 
It is not hard to prove that  $\DX$  is a finite group.

Let $\md$ and $\md '$ be in $\Z^{\gamma}$; we say that they are equivalent 
if  
$\md - \md '\in \LX$.
We denote by $\DX^d$ the set of equivalence classes of multidegrees of total degree $d$.
It is clear that $\DX =\DX^0$ and that
$$
\#\DX =\#\DX^d.
$$

Now back to  $N_X^d$,    the special fiber of  (\ref{neron}); as we said it is a smooth, possibly non connected scheme of pure dimension $g$.

\begin{fact}
\label{number}  Under the above assumptions,
the number of irreducible (i.e. connected) components of $N_X^d$ is equal to $\#\Delta_X$.
\end{fact}
This   is well known; see \cite[8.1.2]{raynaud} (where $\DX$ is  the same as $\ker \beta/\im \alpha$) or   \cite[thm. 9.6.1]{BLR}.
Using the standard notation of N\'eron models theory we have $\DX=\Phi_{N^d_f}$,
i.e. $\DX$ is the ``component group" of $N^d_f$.

\

For every stable curve $X$ and every $d$ we denote by $\PXb$ the degree $d$
compactified Jacobian (or 
degree-$d$   compactified Picard scheme).
$\PXb$ has been constructed in \cite{OS} for a fixed curve, and independently for families
in \cite{simpson} and in \cite{caporaso} (the constructions of \cite{OS} and   \cite{simpson} are here considered    with respect to the canonical    polariz ation); these three constructions give the same scheme by \cite{alex}, see also \cite{LM}.
We mention that another compactified Jacobian is constructed in \cite{est}, whose connection
with the others is under investigation.
 An explicit  description of $\PXb$ will be recalled in \ref{strata}.  We here anticipate the fact that $\PXb$ is a connected projective scheme of pure dimension $g$. As we said in the introduction, several geometric and modular properties of $\PXb$ 
depend on $d$.

\begin{defi}
\label{NDtype} 
Let $X$ be a stable curve and $\PXb$ its degree $d$ compactified Jacobian.
We  say that  $\PXb$ is        of {\em N\'eron type}  if the number of irreducible components of $\PXb$ is equal to the number of irreducible components of $N^d_X$.
\end{defi}

\begin{example}
A curve $X$ is called {\it tree-like} if every node of $X$ lying in two different irreducible components is a separating node.

The compactified Jacobian   of a tree-like curve $X$
is easily seen to be always of N\'eron type. Indeed $\PXb$ is irreducible for every $d$,
on the other hand $\Delta_X=\{0\}$ so that $N_X^d$ is also irreducible.

\end{example}

Let now  $\pi:\pfb\to B$
 be the compactified degree-$d$ Picard scheme of a regular smoothing
  $f:\X\to B$   of $X$, as defined in \ref{neron}. So the fiber of $\pi$ over $b_0$ is
   $\PXb$, and the restriction of $\pi$ over $U=B\smallsetminus \{b_0\}\subset B$ is $\Pic^d_{f_U}$.
 We denote $\pf\to B$ the smooth locus of
$\pi$. By  the N\'eron Mapping Property 
there exists  a canonical $B$-morphism 
 from $\pf$ to the N\'eron model of  $\Pic^d_{f_U}$:
 $$
\chi_f: P^d_f\to \nf
$$
extending the indentity map from the generic fiber of $\pi$ to the generic fiber of  $\nf\to B$.
\begin{prop} 
\label{PN}
In  the above set up,
$\PXb$ is    of   N\'eron type   if and only if the map $\chi_f:\pf \to \nf$ 
 is an isomorphism for every $f:\X\to B$ as above.
 \end{prop}
 The proof, requiring a description of $\pfb$, will be given later, in \ref{proofPN}.
\end{nota}

\begin{nota}{\it Smoothing separating nodes.}
\label{G2}
A stable weighted graph of genus $g$ is a pair  $(\Gamma, w)$,
where $\Gamma$ is a graph and $w:V(\Gamma)\to \Z_{\geq 0}$
a {\it weight function}.
The genus of $(\Gamma, w)$ is the number $g_{(\Gamma, w)}$ defined as follows:
$$
 g_{(\Gamma, w)}=\sum_{v\in V(\Gamma )}w (v)+b_1(\Gamma).
 $$
  A weigthed graph will be called  
{\it stable} if  every $v\in V(\Gamma)$ such that $w(v)=0$ has valency at least $3$.

Let $X$ be a nodal curve of genus $g$, the   weighted dual graph  of $X$ is the weighted graph $(\Gamma_X,w_X)$
such that $\Gamma_X$ is the usual dual graph of $X$
(the vertices of $\Gamma_X$
 are identified with the irreducible components of $X$ and the edges are identified with the nodes of $X$; an edge joins two, possibly equal, vertices if the corresponding node is in the intersection of the corresponding irreducible components),
 and
  $w_X$ is the 
   {\it weight function} on the set of irreducible components of $X$,   $V(\Gamma_X)$,  assigning to a vertex the geometric genus of the corresponding component.
Hence
 $$
 g=\sum_{v\in V(\Gamma_X)}w_X(v)+b_1(\Gamma_X)= g_{(\Gamma, w)}.
 $$
 $X$ is  a stable curve if and only if 
  $(\Gamma_X,w_X)$ is a stable weighted graph.

Now we ask:
What happens to the weighted dual graph of $X$ if we smooth all the separating nodes of $X$?

To answer this question,
we
introduce   a new weighted graph associated to a weighted graph $(\Gamma,w)$.
This will be denoted
$(\Gamma^2,w^2)$ such that 
  $\Gamma^2$ 
is the graph obtained by contracting every separating edge of $\Gamma$ 
to a point. So, $\Gamma^2$ is 2-edge-connected (i.e. free from separating edges); this explains  the notational choice, used elsewhere in the literature.
There is a natural surjective map contracting the separating  edges  of $\Gamma$ 
$$
\sigma: \Gamma\la \Gamma^2
$$
and an induced surjection on the set of vertices
$$
\phi:V(\Gamma)\la V(\Gamma ^2);\  \  \  \  v\mapsto \sigma(v).
$$
The weight function $w^2$ is defined as follows. For every $v^2\in V(\Gamma^2)$
$$
w^2(v^2)=\sum_{v\in \phi^{-1}(v^2)}w(v).
$$

As $\sigma $ does not contract any loop, $b_1(\Gamma)=b_1(\Gamma^2)$ and $g_{(\Gamma, w)}=g_{(\Gamma^2, w^2)}$.
\begin{remark}
\label{X2}
If $(\Gamma,w)$ is the weighted dual graph of a curve $X$,  $(\Gamma^2, w^2)$ is the weighted  dual graph
of any curve obtained by smoothing every separating node of $X$.
We shall usually denote by $X^2$ such a curve. Of course $X$ and  $X^2$ have the same  genus.
\end{remark}
\end{nota}
\begin{nota}{\it $d$-general and weakly $d$-general curves.}
\label{genw}
Let us  recall the definitions of balanced and strictly balanced multidegrees.

\begin{defi}
\label{baldef}
Let $X$ be a  quasistable curve of genus $g\geq 2$ and   $L\in \Pic^dX$. 
Let $\md$ be the multidegree of $L$.
\begin{enumerate}
\item
We say that $L$, or $\md$,  is {\it balanced}  if for any subcurve  (equivalently, for any connected subcurve) 
$Z\subset  Y$  
we have (notation in  \ref{not}(\ref{notZ}))
\begin{equation}
\label{BI}
  \deg_ZL\geq m_Z(d):= \frac{d w_Z}{2g-2}-\frac{\delta_Z}{2},
  \end{equation}
and if $\deg_ZL=1$ if $Z$ is an exceptional component.
 \item
 We say that $L$, or  $\md$,
 is {\it strictly balanced} if it is balanced and if strict inequality holds in  (\ref{BI}) for every  
  $Z\subsetneq  X$ such that $Z\cap Z^c\not\subset \eX$.
\item
We denote  
$$ 
\ov{B_d(X)}=\{\md:|\md|=d  \text{ { balanced on} }  X\}\supset B_d(X)=\{\md:   \text{{  strictly balanced}}\}.
 $$
  \end{enumerate}
\end{defi}

\begin{remark}
\label{in}
Let $X$ be a stable curve. Then every multidegree class in $\DX^d$ has a balanced representative, which is unique if and only if it is strictly balanced; see  \cite[Prop. 4.12]{cner}.
Therefore
$$
\#B_d(X)\leq \#\Delta_X\leq \#\ov{B_d(X)}.
$$
\end{remark}

The terminology ``strictly balanced"    is not to be confused with   ``stably balanced"
(used elsewhere and unnecessary here).
The two coincide for stable curves,
but, in general, a stably balanced line bundle is strictly balanced, but the converse may fail.

Let us explain the difference.
The compactified Picard scheme of $X$, $\PXb$
is   a GIT-quotient of a certain scheme by a certain group $G$.
Strictly balanced line bundles correspond to the  GIT-semistable orbits
that are closed in the GIT-semistable locus. Stably balanced line bundles correspond to GIT-stable points and
balanced line bundles correspond to GIT-semistable points.
As 
 every point in $\PXb$ parametrizes a unique closed orbit,
 strictly balanced line bundles of degree $d$ on quasistable curves
 of $X$
 are bijectively parametrized by 
$\PXb$. See Fact~\ref{strata} below.

Recall from \cite[4.13]{cner} that a stable curve $X$ is called $d$-general if 
$B_d(X)=\ov{B_d(X)}$.
Equivalently, $X$ is $d$-general if the inequalities in Remark~\ref{in} are all equalities.

\begin{remark}
\label{dgen}
The following facts are  well known:
\begin{enumerate}
\item
The set of $d$-general stable curves is a nonempty open subset of $\Mgb$.
\item
 $(d-g+1,2g-1)=1$  if and only if every stable curve of genus $g$ is $d$-general.
\item
The property of being $d$-general depends only on the weighted dual graph
(obvious).
\end{enumerate}
\end{remark}
\begin{defi}
\label{wdef}
Let $X$ be a stable curve. We will say that $X$, or its weighted dual graph $(\Gamma_X, w_X)$, is
$d${\it -general} if $B_d(X)=\ov{B_d(X)}$.

We will say that $X$ is {\it weakly d-general} if $(\Gamma_X^2, w_X^2)$ is $d$-general.
\end{defi}

\begin{example}
If $\sep =\emptyset$ then $X$ is $d$-general if and only if it is weakly general.

If $X$ is tree-like, then  $(\Gamma_X^2, w_X^2)$ has only one vertex, hence
 it is $d$-general for every $d$. Therefore tree-like curves are weakly $d$-general for every $d$.
\end{example}
\end{nota}
\section{Irreducible components of compactified Jacobians}
\label{sec2}
\begin{nota}{\it Compactified degree-$d$ Jacobians.}
Let us describe the compactified Jacobian $\PXb$ for any degree $d$. 
We use the set up of  \cite{caporaso} and \cite{cner}; 
in these papers there is the assumption $g\geq 3$, but by
\cite{OS},  \cite{simpson} and \cite{alex} we can extend our results to $g\geq 2$.
A synthetic account of the modular properties of the compactified Jacobian for a curve
or for a family can be found in \cite[3.8 and 5.10]{CE}.
\begin{fact}
\label{strata}
Let $X$ be a stable curve of genus $g\geq 2$.
Then $\PXb$ is a connected, reduced, projective scheme of pure dimension $g$, admitting a canonical decomposition (notation in \ref{not}(\ref{notS}))
$$
\PXb=\coprod_{\stackrel{S\subset \sing}{\md \in B_d(\hXS)}} P^{\md}_S
$$
such that for every $S\subset \sing $ and $\md\in  B_d(\hXS)$ there is a natural isomorphism
$$ P^{\md}_S\cong \Pic^{\mdn}\nXS
$$
where $\mdn$ denotes the multidegree on  $\nXS\subset \hXS$
defined by restricting
 $\md$.

Let ${\rm{i}}(\PXb)$ be the number of irreducible components of $\PXb$; then
\begin{equation}
\label{ineq}
B_d(X)\leq {\rm{i}}(\PXb)\leq \#\Delta_X.
\end{equation}
\end{fact}

\begin{cor}
\label{bij} Let $X$ be a stable curve.
\begin{enumerate}
\item
The decomposition  of $\PXb$ in  irreducible components is
$$
\PXb =\bigcup_{(S,\md)\in I^d_X}\ov{P^{\md}_S},\  \  \text{ where }  \  \  I^d_X:=\{(S,\md): \      S\subset \sep,\   \md \in B_d(\hXS)\}.
$$
\item
\label{genner}
Suppose that $X$ is $d$-general; then    $\PX$  is of N\'eron type, and
for every nonempty $S\subset \sep$ we have $B_d(\hXS)=\emptyset$
\end{enumerate}
\end{cor}
\begin{proof}
From   Fact~\ref{strata} we have that 
the irreducible components of $\PXb$
are the closures of subsets $P^{\md}_S\cong \Pic^{\mdn}\nXS$ where $S$ is such that $\dim\Pic^{\mdn}\nXS=g$.
Now, it is clear that 
$$
\dim \Pic^{\mdn}\nXS=\dim J(\nXS) =g \ \   \text{ if and only if}\  \  \  \  S\subset \sep.
$$
Therefore the irreducible components of $\PXb$ correspond bijectively
to pairs $(S,\md)$ with $S\subset \sep$ and $\md \in B_d(\hXS)$.

Now part (\ref{genner}). 
It is clear that the set $I^d_X$ contains  a subset identifiable with $B_d(X)$,
namely the subset $\{(\emptyset ,\md): \        \md \in B_d(X)\}$.
If $X$ is $d$-general then $ \#B_d(X)=\#\Delta_X$, hence by (\ref{ineq})
we must have that $I^d_X$ contains no pairs other than those of type $(\emptyset ,\md)$.
This concludes the proof.
\end{proof}

\begin{lemma}
\label{mu}
Let $X$ be a stable curve of genus $g$ and 
let $\mu\in \Delta_X^d$ be a multidegree class.
Then there exists a unique $S(\mu)\subset \sing$ and a unique $\md(\mu)\in B_d(\widehat{X_{S(\mu)}})$ such that for every $\md \in \ov{B_d(X)}$ with
$[\md]=\mu$ 
  the following properties hold.
\begin{enumerate}
\item
There is a canonical surjection
$$
\Pic^{\md}X \la P^{\md(\mu)}_{S(\mu)}\cong \Pic^{\md(\mu)^{\nu}}X^{\nu}_{S(\mu)}.
$$
\item
\label{Smu} We have
$$
S(\mu)=\bigcup_{Z\subset X:d_Z=m_Z(d)} Z\cap Z^c.
$$
\end{enumerate}
\end{lemma}
\begin{proof}
The proof is standard. 
Let us sketch it using the combinatorial results   \cite[Lemma 5.1 and Lemma 6.1]{caporaso}.
The terminology used in that paper  differs from ours as follows:
what we here call a ``strictly balanced multidegree $\md$ on a quasistable curve $X$"
is there called an ``extremal pair $(X,\md)$"; cf. subsection 5.2 p. 631.

So, 
the pair    $(\widehat{X_{S(\mu)}}, \md(\mu))$
is the ``extremal pair" associated to $\mu$. This means the following.
For every balanced line bundle $L$ on $X$ such that $[\mdeg L]=\mu$
the    point in $\PXb$ associated to $L$ parametrizes a  line bundle
$\hL\in \Pic^{\md} \widehat{X_{S(\mu)}}$, and the restriction
of $\hL$ to $\XS$ is uniquely determined by $L$. 
Conversely every line bundle in $\Pic^{\md(\mu)^{\nu}}X^{\nu}_{S(\mu)}$
is obtained in this way.

More precisely, recall that $\PXb$ is a GIT quotient, let us denote it by $\PXb=V_X/G$
so that $V_X$ is made of GIT-semistable points, and
$G$ is a reductive group. Let $O_G(L)\subset V_X$ be the orbit of $L$.
Then the semistable closure of $O_G(L)$ contains a unique closed orbit
$O_G(\hL)$ as above.
Moreover for every $\md'\in B_d(X)$  having class $\mu$ there exists $L'\in \Pic^{\md'}X$
such that the above $O_G(\hL)$ lies in the closure of $O_G(L')$.
Hence the maps $\Pic^{\md}X\to \PXb$ and $\Pic^{\md'}X\to \PXb$ have the same image.
Parenthetically,  this ensures separatedness of the relative compactified Picard scheme; see \cite[Sec. 3]{cner}.

Using the notation  of \ref{strata}, we have that for every
balanced $\md$   of class $\mu$  the  canonical map  $\Pic^{\md}X\to \PXb$ has
 image 
$P^{\md(\mu)}_{S(\mu)}$, so that the first part is proved. 

Now  (\ref{Smu}).
The previously mentioned Lemma 5.1 
 implies that
  for every $\md\in  \ov{B_d(X)}$ and every $Z$ 
  such that $\md_Z-m_Z(d)$ we have $Z\cap Z^c\subset S(\mu)$.
By the above Lemma 6.1 each $n\in S(\mu)$ is obtained in this way
  \end{proof}

\begin{prop}
\label{prop}
Let $X$ be a stable curve of genus $g$. $\PXb$ is of N\'eron type if and only if for every $\md \in   \ov{B_d(X)}$ and every connected  $Z\subset X$ such that $d_Z=m_Z(d)$ we have $Z\cap Z^c\subset \sep$.
\end{prop}
\begin{proof}
We begin by observing that, with the notation of Corollary~\ref{bij} and Lemma~\ref{mu},
we have
$$
I^d_X=\{(\md(\mu),S(\mu)),\     \forall \mu\in \DX^d\  \text{ such that}\   \dim P^{\md(\mu)}_{S(\mu)}=g\}.
$$
Indeed, by Fact~\ref{strata} the set on the right is clearly included in $I^d_X$.
On the other hand let $(S,\md)\in I^d_X$.
To show that there exists $\mu\in \DX^d$ such that $\md=\md(\mu)$ we can   assume that $S\neq \emptyset$ (otherwise it is obvious).
So, $\md$ is a strictly balanced multidegree of total degree $d$ on $\hXS$.
Let $n\in S$; by Corollary~\ref{bij} $n$ is a separating node of $X$; let $X=Z\cup Z^c$
with $Z\cap Z^c=\{n\}$. Then $Z$ and $Z^c$ can also be viewed as subcurves of
$\hXS$, where they do not intersect since the node $n$ is replaced by an exceptional component $E$. Now, $\md_E=1$, therefore $\md_Z=m_Z(d)$ and  $\md_{Z^c}=m_{Z^c}(d)$.
Let $C_Z\subset Z\subset X$ be the irreducible component intersecting $Z^c$
(so that  $C_Z\subset \hXS$ intersects $E$).
Let $ {\md^X}$ be the multidegree on $X$ defined  
as follows: for every irreducible component $C\subset X$
\begin{displaymath}
\md^X_C=\left\{\begin{array}{ll}
\md_C+1 &\text{ if } C=C_Z,\\
\md_C &\text{ otherwise.}\\
\end{array}\right.
\end{displaymath}
As $\md$ is balanced on $\hXS$ one easily cheks that $\md^X$ is balanced on $X$.
Note that  $\md^X$ is not strictly balanced, since  $\md^X_{Z^c}=m_{Z^c}(d)$.
By iterating the above procedure for every node in $S$ we arrive at a balanced
multidegree on $X$ whose class we denote by $\mu\in \DX^d$.
By Lemma~\ref{mu} we have that $\md =\md(\mu)$.

 Suppose that  $\PXb$ is of N\'eron type. By the previous discussion 
there is a natural bijection between $I^d_X$ and $\Delta_X^d$, mapping
  $\mu\in \Delta_X^d$
 to $(S(\mu),\md(\mu))$.
By Corollary~\ref{bij} we have  $S(\mu)\subset \sep$. Hence
 for every multidegree $\md\in   \ov{B_d(X)}$
such that $[\md]=\mu$ we have that condition (\ref{Smu}) of that lemma holds.
In particular every $Z$ as in our statement is such that $Z\cap Z^c\subset S(\mu)\subset \sep$.

Conversely, if $\PX$ is not of N\'eron type there is a class $\mu\in \Delta_X^d$ such that 
$$
g> \dim P^{\md(\mu)}_{S(\mu)}=\dim  J(X^{\nu}_{S(\mu)}).
$$
But then $S(\mu)$ must contain some non separating node of $X$. Therefore, by Lemma~\ref{mu}(\ref{Smu}),  there exists a connected subcurve $Z\subset X$ such that $d_Z=m_Z(d)$ and such that $Z\cap Z^c$ contains some non separating node. 
\end{proof}
\end{nota}
\begin{nota}{\it Proof of Proposition~\ref{PN}.}
\label{proofPN}
We generalize the proof of \cite[Thm. 6.1]{cner}.
Let $f:\X \to B$ be a regular smoothing
  of $X$ as defined in \ref{neron} and $\pi:\pfb\to B$ be the compactified degree-$d$ Picard scheme.
Its smooth locus $\pf\to B$ is such that  its fiber  over $b_0$, denoted $ P^d_X$, satisfies
\begin{equation}
\label{PXd}
   P^d_X=\coprod_{(S,\md)\in I^d_X}  P^{\md}_S
\end{equation}
   (notation in \ref{bij}) where each $P^{\md}_S$ is irreducible of dimension $g$.
   If the morphism $\ph_f:\pf \to N_f^d$ is an isomorphism,
   then $P^d_X$ has as many irreducible components as $N_X^d$, hence the same holds for $\PXb$. So $\PXb$ is of N\'eron type.
   
Conversely, if $\PXb$ is of N\'eron type, then $P^d_X$ has an irreducible component for
   every $\mu\in \DX^d$ so that   (\ref{PXd}) takes the form
   $$
   P^d_X=\coprod_{\mu \in \DX^d}  P^{\md(\mu)}_{S(\mu)}.
   $$
Let us construct the inverse of $\chi_f$.
   We pick a balanced representative $\md^{\mu}$ for every multidegree class $\mu\in \DX^d$
   (it exists by Remark~\ref{in}).
   By  \cite[Lemma. 3.10]{cner} we have 
   $$
   N^d_f\cong \frac{\coprod _{\mu \in \DX^d} \Pic^{\md^{\mu}}_f}{\sim _U}
   $$
   where $\sim _U$ denotes the gluing of the Picard schemes $\Pic^{\md^{\mu}}_{f}$
   along their restrictions over $U$ (as $\Pic^{\md^{\mu}}_{f_U} =\Pic^d_{f_U}$ for every $\mu$).
Now,  the Picard scheme $\Pic^{\md^{\mu}}_f$ is endowed with  
a Poincar\'e bundle, which is a relatively balanced   line bundle on
$\X\times _B\Pic^{\md^{\mu}}_f$. By the modular property of $\pfb$ the  Poincar\'e bundle induces a   canonical $B$-morphism 
$$
\psi_f^{\mu}:\Pic^{\md^{\mu}}_f \la P^{\md(\mu)}_{S(\mu)}\subset \pfb.
$$
As $\mu$ varies, the restrictions of these morphisms over $U$
all coincide with the identity map $\Pic^d_{f_U}\to \Pic^d_{f_U}\subset \pfb$.
   Therefore  the $\psi_f^{\mu}$   can be glued together to a morphism
   $$
   \psi_f:N^d_f\la \pf\subset \pfb.
   $$
   It is clear that $\psi_f$ is the inverse of $\chi_f$. Proposition~\ref{PN} is proved. \hfill$\qed$
\end{nota}
\begin{nota}{\it The main result.}
From Proposition~\ref{prop} we derive the following.
\begin{cor}
\label{nosep}
Let $X$ be a stable curve free from separating nodes.
Then $\PXb$ is of N\'eron type if and only if $X$ is $d$-general.
\end{cor}
\begin{proof}
Suppose  $\sep=\emptyset$.
Then if there exists $\md\in \ov{B_d(X)}\smallsetminus B_d(X)$, such a $\md$ will not satisfy the condition of \ref{prop}. Hence if $X$ is not $d$-general,  $\PXb$ is not of N\'eron type.
The converse follows from Corollary~\ref{bij}(\ref{genner}).
\end{proof}

We are ready to prove our main result.
\begin{thm}
\label{main}
Let $X$ be a stable curve.
Then $\PXb$ is of N\'eron type if and only if $X$ is weakly $d$-general.
\end{thm}
\begin{proof}
Observe that if $X$ is free from separating nodes we are done by Corollary~\ref{nosep}.
Let $(\Gamma ,w)$ be the
  weigthed graph of $X$ and consider the weighted graph $(\Gamma^2,w^2)$
  defined in subsection~\ref{G2}.
We denote by $X^2$ a stable curve
whose weighted graph is  $(\Gamma^2,w^2)$; see Remark~\ref{X2}.

  Recall that we denote by $\sigma:\Gamma \to \Gamma^2$ the   contraction map
  and by 
  $$
  \phi:V(\Gamma)\to V(\Gamma^2);\  \  \  v\mapsto \sigma(v)
  $$ 
  the induced map on the vertices, i.e.  on the irreducible components.
  The subcurves of $X$ naturally correspond to the so-called ``induced" subgraphs  of $\Gamma$,
  i.e. those subgraphs $\Gamma'$  such that if two vertices $v,w$ of $\Gamma$ are in $\Gamma'$,
  then    every edge of $\Gamma$   joining $v$ with $w$ lies in
$\Gamma'$.
Similarly, the induced subgraphs of $\Gamma^2$ correspond to subcurves   of $X^2$.
If $Z^2$ is a subcurve of $X^2$,
and $\Gamma_{Z^2}\subset \Gamma_{X^2}$ its corresponding subgraph,
we denote by $Z\subset X$ the subcurve associated to
$\sigma^{-1}(\Gamma_{Z^2})$ (it is obvious that the subgraph  $\sigma^{-1}(\Gamma_{Z^2})$ is   induced if so is $\Gamma_{Z^2}$);
we refer to $Z$ as the ``preimage" of $Z^2$. 
Of course  $\sigma(\Gamma_Z)=\Gamma_{Z^2}$.

For any  $Z\subset X$ which is the preimage of a subcurve $Z^2\subset X^2$ we have  
\begin{equation}
\label{Zsep}
Z\cap \sep \subset Z_{\rm{sep}} 
\end{equation}
or, equivalently, $Z\cap Z^c\cap \sep =\emptyset$.
Conversely, every $Z\subset X$ satisfying (\ref{Zsep}) is the preimage of some $Z^2\subset X^2$.

Hence  $Z^2$ can be viewed as a smoothing of $Z$ at its separating nodes
that are also separating nodes of $X$, i.e. at $Z_{\text{sep}}\cap \sep$.
Thus, for every $Z^2$ with preimage $Z$ we have $g_Z=g_{Z^2}$ and  $\delta_Z=\delta_{Z^2}$;  hence for every $d\in \Z$
$$
m_{Z^2}(d)=m_{Z}(d).
$$
A multidegree can, and will, be viewed as an integer valued function on the vertices.
We claim that there is a surjection
$$
\alpha:\ov{B_d(X)}\la \ov{B_d(X^2)}
$$
such that 
$$
\forall v^2\in V(\Gamma^2),\  \  \  \  \  \  \alpha(\md)(v^2)=\sum_{v\in \phi^{-1}(v^2)}\md (v).
$$

If $\md$ is balanced, so is $\alpha(\md)$, 
indeed for every subcurve $Z^2\subset X^2$ we have $\alpha(\md)_{Z^2}=\md_Z$
 where $Z\subset X$ is the preimage of $Z^2$;
  by what we said above  the inequality (\ref{BI}) is satisfied on $Z^2$ if (and only if)  it is satisfied on $Z$.

Let us now show that $\alpha$ is surjective.
Let $\md^2$ be a balanced multidegree on $X^2$; recall that $X^2$ can be chosen to be a smoothing of $X$ at all of its separating nodes. In other words there exists a family of curves
$X_t$, all having $(\Gamma^2,w^2)$ as weighted graph, specializing to $X$.
But then there also exists a family of balanced line bundles $L_t$ on $X_t$
specializing to a line bundle of degree $\md$ on $X$ (this follows from the construction of the universal compactified Picard scheme $\Pdgb\to \Mgb$,
see \cite[subsection 5.2]{cner} for an overview).
By the definition of $\alpha$, 
it is clear that the multidegree $\mdeg L_t$ is such that $\alpha(\mdeg L_t)=\md^2$.

Now assume that $\PXb$ is of N\'eron type.
Our goal is to prove that $X^2$ is $d$-general. By contradiction, let $Z^2\subset X^2$
be a connected subcurve such that for some $\md^2 \in\ov{B_d(X^2)}$ we have $\md^2_{Z^2}=m_{Z^2}(d)$.
Let $Z$ be the preimage of $Z^2$, and let $\md\in \ov{B_d(X)}$ be such that $\alpha(\md)=\md^2$. Then 
$$
\md_Z=\md^2_{Z^2}=m_{Z^2}(d)=m_{Z}(d).
$$
By Proposition~\ref{prop} we obtain that $Z\cap Z^c\subset \sep$. Now this is not possible, 
indeed,
as $Z$ is the preimage of $Z^2$,    (\ref{Zsep}) holds, and therefore 
$Z\cap Z^c\not\subset \sep.$

Conversely, let $X$ be a weakly $d$-general curve.
To show that $\PX$ is of N\'eron type we use
again Proposition~\ref{prop}, according to which it suffices to show that for every $\md\in  \ov{B_d(X)}$ and for every $Z\subset X$ such that $Z\cap Z^c\not\subset \sep$ we have $\md_Z>m_Z(d)$.

By contradiction. Let $Z$ be a connected subcurve such that $Z\cap Z^c\not\subset  \sep$,
and $\md_Z=m_Z(d)$ for some balanced multidegree $\md$ on $X$. 
We choose $Z$ maximal with respect to this property, so that we also have 
\begin{equation}
\label{Zmax}
Z\cap Z^c\cap \sep =\emptyset.
\end{equation}
Indeed, if $Z\cap Z^c$ contains a separating node $n\in \sep$, there is a connected component $Z'\subset Z^c$ such that $Z\cap Z'=\{n\}$; by replacing $Z$ with $W:=Z\cup Z'$
we have a connected curve $W$ containing $Z$ and such that $W\cap W^c\not\subset\sep$;
a trivial computation shows    that $\md_W=m_W(d)$, contradicting the maximality of $Z$.

By (\ref{Zmax})
we have that $Z\cap  \sep$ is all contained in $Z_{\rm{sep}}$ therefore, as observed in (\ref{Zsep}),
the curve $Z$ is the preimage of a subcurve $Z^2\subset X^2$
(notation as in the previous part of the proof). 
Now let 
$\md^2=\alpha(\md)$; so $\md^2\in \ov{B_d(X^2)}=B_d(X^2)$ by hypothesis.
We have
$$
\md^2_{Z^2}=\md _Z=m_Z(d)=m_{Z^2}(d).
$$
This contradicts the fact that $\md^2$ is strictly balanced on $X^2$. 
\end{proof}

\begin{cor}
\label{g-1D}
Let $X$ be a stable curve.
If $d=g-1$ then $\PXb$  is of N\'eron type if and only if $X$ is a tree-like curve.
\end{cor}
\begin{proof}
Let $d=g-1$. Then   $X$ is $d$-general if and only if $X$ is irreducible
(see \cite[Rk. 2.3]{melo}).
Hence $X$ is weakly $d$-general if and only if $X$ is tree-like.
\end{proof}
\end{nota}
\begin{nota}
\label{locus}
{\it The locus of weakly $d$-general curves in $\Mgb$.}
The locus of $d$-general curves in $\Mgb$ has been studied in details in \cite{melo} (also in \cite{CE} if $d=1$ for applications to Abel maps).
A stable curve $X$ which is not $d$-general is called {\it d-special}.
The locus of $d$-special curves is a closed subscheme denoted $\Sigma_g^d\subset \Mgb$.
By \cite[Lemma 2.10]{melo}, $\Sigma^d_g$ is the closure of the locus of $d$-special curves made of two smoth components.
Curves made of two smoth components are called {\it vine curves}. We are going to exhibit a precise description of the complement of the locus of 
weakly $d$-general curves in $\Mgb$, i.e. of the locus of curves whose degree-$d$ compactified Jacobian is not of N\'eron type.
We denote it by $\Ddg\subset \Mgb$. i.e.
$$
\Ddg:=\{X\in \Mgb: \  \PXb \text{ not of N\'eron type}\}.
$$
 In the following statement by $\codim\Ddg$ we mean the codimension of an irreducible component of maximal dimension.
 \begin{prop}
 \label{D}
  $\Ddg$ is  the closure of the locus of $d$-special vine curves with at least $2$ nodes.
Moreover
\begin{displaymath}
\codim \Ddg =\left\{\begin{array}{ll}
+\infty    \   (\text{i.e. } \Ddg=\emptyset) & \text{ if } (d-g+1,2g-2)=1,\\
3&\text{ if } (d-g+1,2g-2)=2  \text{ and }g \text{ is even, }\\
2 &\text{ otherwise.}\\
\end{array}\right.
\end{displaymath}

\end{prop}
\begin{proof}
By Theorem~\ref{main},
we have that $X\in \Ddg$ if and only if $X$ is not weakly $d$-general, if and only if $X^2$ is not $d$-general (where $X^2$ is as in \ref{X2}).
This is equivalent to the fact that there exists $\md \in B_d(X^2)$ and a subcurve $Z\subset X^2$ such that $\md_{Z}=m_Z(d)$; as $X^2$ has no separating nodes, for every subcurve $Z\subset X^2$ we have
$\delta_Z\geq 2$.
This observation added to the proof of \cite[Lemma 2.10]{melo} gives that $X^2$
(and every curve with the same weighted graph) lies in the closure of the locus of $d$-special vine curves with at least two nodes.
Therefore the same holds for $X$, since $X$ is a specialization of curve with the same 
weighted graph as $X^2$.

Conversely, let $X$ be in the closure of the locus of $d$-special vine curves with at least two nodes.
Then $X^2$ is also in this closure, as such vine curves are obviously free from separating nodes.
By \cite[Lemma 2.10]{melo} the curve $X^2$ is $d$-special, hence $X$ is not weakly $d$-general.

Let us turn to $\codim \Ddg$.  
The fact that if $(d-g+1,2g-2)=1$ then  
$\Ddg$ 
  is empty   is well known (\cite{cner}).
  Conversely, assume $\Ddg=\emptyset$. By the previous part, the locus of $d$-special vine curves with at least two nodes is also empty.
  Now the proof of the numerical Lemma 6.3 in \cite{caporaso} shows that this implies that  $(d-g+1,2g-2)=1$. Indeed
  the proof  of that Lemma shows that if there are no $d$-special vine curves with two or three nodes then $(d-g+1,2g-2)=1$. 
  
  Next, recall that the locus of vine curves with $\delta$ nodes has pure codimension $\delta$,
  and notice that the sublocus of $d$-special curves is a union of irreducible components.

Now, again by the proof of the above Lemma 6.3,
if $(d-g+1,2g-2)\neq 1$  and if there are no $d$-special vine curves with two nodes,
  then $(d-g+1,2g-2)=2$, $g$ is even and every vine curve with three nodes, having one component of genus $g/2-1$, is $d$-special.
  This completes the proof of the Proposition.
\end{proof}
A precise description of the locus of $d$-special vine curves is given in \cite[Prop.  2.13]{melo}.
Her result combined with the previous proposition yields a  more precise description of the locus of stable curves
whose compactified degree-$d$ Jacobian is of N\'eron type, for every fixed $d$.
\end{nota}

\end{document}